\documentclass[11pt]{article}
\usepackage{amsmath}
\usepackage{verbatim}
\usepackage{geometry}                
\usepackage{graphicx}
\usepackage{amssymb}
\usepackage{epstopdf}
\newtheorem{thm}{Theorem}[section]
\newtheorem{cor}[thm]{Corollary}
\newtheorem{lem}[thm]{Lemma}

\title{Two Kinds of Division Polynomials
For Twisted Edwards Curves}
\author{Richard Moloney  
and
Gary McGuire\thanks{Research supported by Claude Shannon Institute,
Science Foundation Ireland Grant 06/MI/006, and Grant 07/RFP/MATF846, 
and the Irish Research Council
for Science, Engineering and Technology}\\
School of Mathematical Sciences\\
University College Dublin\\
Ireland\\
Email: \texttt{richard.moloney@ucd.ie}, \texttt{gary.mcguire@ucd.ie}}
\begin{document}

\maketitle

\begin{abstract}
This paper presents two kinds of division polynomials for twisted Edwards curves. Their chief property is that they characterise the $n$-torsion points of a given twisted Edwards curve. 
We also present results concerning the coefficients of these polynomials, which may aid computation.
\end{abstract}

\section{Introduction}


The famous last entry in the diary of Gauss concerns the curve with equation
\begin{equation}\label{gauss1}
x^2+y^2+x^2y^2=1
\end{equation}
and its rational points over $\mathbb{F}_p$. This curve is related to the elliptic curve $y^2=4x^3-4x$.

The idea of \emph{division polynomials} on a curve with a group law on its points, is that we try to
write down a formula for $[n]P$ in terms of the coordinates of $P$,
where $[n]P$ denotes $P$ added to itself $n$ times under the group law.
In this paper we shall give two distinct solutions to this problem, in the general context of 
twisted Edwards curves, of which (\ref{gauss1}) is a special case.

Edwards  \cite{Edwards}, generalising (\ref{gauss1}), introduced an addition law on the curves $x^2+y^2 = c^2(1+x^2y^2)$ 
for $c \in k$, where $k$ is a  field of characteristic not equal to 2. He showed that every elliptic curve over $k$ is birationally equivalent (over some extension of $k$) to a curve of this form. 
 
In \cite{BL}, Bernstein and Lange generalised this addition law to the curves 
$x^2+y^2 = 1+ dx^2y^2$ for $d \in k \setminus \{0,1\}$. 
More generally, they consider $x^2 + y^2 = c^2(1+dx^2y^2)$, however, any such curve is isomorphic to one of the form  $x^2+y^2 = 1+ d^{\prime}x^2y^2$ for some $d^{\prime} \in k$, 
so we will assume $c = 1$. These curves are referred to as Edwards curves. Bernstein and Lange showed that if $k$ is finite, a large class of elliptic curves over $k$ (all those which have a point of order 4) can be represented in Edwards form.
The case $d=-1$ gives the curve (\ref{gauss1}) considered by Gauss.

In \cite{Twisted}, Bernstein et al. introduced the twisted Edwards curves $ax^2+y^2 = 1+dx^2y^2$ (where $a, \ d \in k$ are distinct and non-zero) and showed that every elliptic curve with a representation in Montgomery form is birationally equivalent to a twisted Edwards curve.
Obviously, the case $a=1$ of a twisted Edwards curve is an Edwards curve.

In this paper we describe a sequence of rational functions, and consequently a sequence of
polynomials, defined on the function field of a twisted Edwards curve which are analogous to the division polynomials for elliptic curves in Weierstrass form. 
 In particular, these 
polynomials characterise the $n$-torsion points of the 
twisted Edwards curve for  a positive integer $n$ (see Corollary \ref{divformcor}
and Corollary \ref{coronevar}).
These twisted Edwards division polynomials are polynomials in $y$ with coefficients in $\mathbb{Z}[a,d]$, and have degree in $y$ less than $n^2/2$.

This paper is laid out as follows.
In Section 2 we recall division polynomials for elliptic curves in Weierstrass form.
Section 3 recalls the basic properties of twisted Edwards curves.
In Section 4, on the function field of an Edwards curve, Theorem \ref{uniq1}
proves a uniqueness form for elements of the function field of an Edwards curve,
analagous to the known result that elements of the function field of a
Weierstrass curve can be written uniquely in the form $p(x)+yq(x)$.
Our division polynomials (actually rational functions) are presented in this unique form.
Section \ref{gauss} compares our results to those of Gauss for the curve  (\ref{gauss1}).
In Section \ref{divpols} we isolate the important part of the Edwards division rational functions, 
which are polynomials that could be called Edwards division polynomials.
Furthermore, we show in Section \ref{further}
 that the coefficients of a given twisted Edwards division polynomial exhibit a certain symmetry, which may reduce the amount of computation necessary for finding that polynomial.
In Section \ref{alt}, we derive a different set of polynomials which also display some properties we require from division polynomials. These have a different character to the first set, since the $n$th polynomial is defined by a recursion on the $n-1$th and $n-2$th polynomials, as opposed to polynomials of index $\sim \frac{n}{2}$.

\section{Division polynomials for Weierstrass Curves} 

We recall the division polynomials for Weierstrass curves here.

First we recall the definition of the function field of an (affine) algebraic variety. 
If $V/k$ is a variety in affine $n$-space, 
$I(V)$ denotes the ideal generated by the polynomials in $k[x_1,\dots, x_n]$
that vanish on $V$.
The affine coordinate ring of $V$ is the integral domain 
\[
k[V]:=k[x_1,\dots, x_n]/I(V).
\] 
The function field of $V$ over $k$, denoted by $k(V)$, is defined to be the quotient field of $k[V]$.

For example, if $W$ is an elliptic curve with Weierstrass equation $v^2 = u^3+Au+B$, the function field of $W$, $k(W)$, is the quotient field of $k[u, v]/(v^2-u^3-Au-B)$.

We use $(u,v)$ as the coordinates for a curve in Weierstrass form and reserve $(x,y)$ for (twisted) Edwards curves.

If $char(k) \neq 2$ or 3, given an elliptic curve over $k$ in short Weierstrass form 
\[
W: v^2 = u^3+Au+B
\] 
with identity $\mathcal{O}$ , the division polynomials $\Psi_{n}$ are polynomials defined on the function field of $W$ for each $n \in \mathbb{N}$ by the following recursion:
\begin{align*}\Psi_{0}(u, v) &= 0\\
\Psi_{1}(u, v) &= 1\\
\Psi_{2}(u, v) &= 2v\\
\Psi_{3}(u, v) &= 3u^4+6Au^2+12Bu - A^2\\
\Psi_{4}(u, v) &= 4v(u^6+5Au^4+20Bu^3-5A^2u^2-4ABu-A^3-8B^2)\\
\Psi_{2m+1}(u, v) &= \Psi_{m+2}(u,v)\Psi_{m}^3(u,v) - \Psi_{m-1}(u,v)\Psi_{m+1}^3(u,v) \text{ for } m \geq 2\\
\Psi_{2m}(u,v) &= {\Psi_{m}(u,v)\over\Psi_{2}(u,v)}\left(\Psi_{m+2}(u,v)\Psi_{m-1}^2(u,v) -\Psi_{m-2}(u,v)\Psi_{m+1}^2(u,v)\right)\quad \text{ for }m \geq 3.\end{align*}

The $\Psi_{n}$ are polynomials in $u$ and $v$ with coefficients in $\mathbb{Z}[A,B]$. The principal properties of the division polynomials are that $\Psi_{n}(u,v) = 0$ precisely when $(u,v)$ is an $n$-torsion point of $W$ (i.e. $[n](u,v) = \mathcal{O}$), and that the multiplication-by-$n$ map $[n]:W \rightarrow W$ is characterised by the division polynomials as 
\[ [n](u,v) = \left({{u\Psi_{n}^2(u,v) - \Psi_{n-1}(u, v)\Psi_{n+1}(u, v)}\over{\Psi_{n}^2(u, v)}}, {{\Psi_{2n}(u, v)}\over{2\Psi_{n}^4(u, v)}}\right)\] 
(see e.g. \cite{Wash}, Chapters 3 , 9, \cite{Silv}, Chapter 3). 
If $n$ is odd then  $\Psi_{n} \in \mathbb{Z}[u,A,B]$, and $\Psi_n$
has degree  $(n^2-1)/2$ in $u$.  If $n$ is even then
$\Psi_{n} \in v\mathbb{Z}[u,A,B]$ with degree  $(n^2-4)/2$ in $u$.
In this paper we prove analagous results for twisted Edwards curves.

\section{Twisted Edwards Curves}\label{tecsection}

Let $k$ be a field with characteristic $\neq 2$ or $3$. 
Let $K$ be an extension field of $k$.
Let $E(K)$ be the twisted Edwards curve over $K$
with coefficients $a$ and $d$, where $a$ and $d$ are distinct and non-zero:  
\[E(K): ax^2+y^2 = 1+dx^2y^2.
\]
Points on $E(K)$ may be added by the rule 
\[(x_1, y_1)+(x_2, y_2) = \left({{x_1y_2+x_2y_1}\over{1+dx_1x_2y_1y_2}} , {{y_1y_2-ax_1x_2}\over{1-dx_1x_2y_1y_2}}\right)\]
and under this operation, the points on $E(K)$ 
form an abelian group.
The identity is $(0,1)$, and the additive inverse of a point $(x,y)$ is $(-x,y)$.
The projective closure of $E$ has singularities at  $(1:0:0)$ and $(0:1:0)$.

The twisted Edwards curve
$E(K)$ is birationally equivalent to the Weierstrass-form elliptic curve
\[W(K):  v^2 = u^3 - {{(a^2+14ad+d^2)}\over{48}}u - {{(a^3-33a^2d-33ad^2+d^3)}\over{864}}\] under the transformation 
\[
u := {{(5a-d)+(a-5d)y}\over{12(1-y)}} \ , \  v := {{(a-d)(1+y)}\over{4x(1-y)}} \quad \text{if } x(1-y)\neq 0, \]
otherwise
\begin{align*}(x,y) = (0,1) &\Rightarrow (u, v) = \mathcal{O}\\
(x,y) = (0,-1) &\Rightarrow (u, v) = \left({{a+d}\over{6}}, 0\right).\end{align*}

The inverse transformation is given by
\[x = {{6u-(a+d)}\over{6v}}, \  y = {{12u+d-5a}\over{12u+a-5d}} \quad \text{if } v(12u+a-5d)\neq 0  \]
and
\begin{align*}
(u,v) = \mathcal{O} &\Rightarrow (x,y) = (0,1)\\
(u,v) = \left({{a+d}\over{6}}, 0\right) &\Rightarrow (x,y) = (0,-1).
\end{align*}

There are 4 points on $W(\overline{k})$ that are not mapped to any point on the twisted Edwards curve. These are $(u, v) = \left({{5d-a}\over{12}} , \pm{{s(d-a)}\over{4}}\right)$ and $(u, v) =\left({{-(a+d) \pm 6t}\over{12}},0\right) $ where $s, t \in \bar{k}$ such that $s^2 = d, t^2 = ad$. We note that $\left({{-(a+d) \pm 6t}\over{12}},0\right)$ are points of order 2 on $W$, and $\left({{5d-a}\over{12}} , \pm{{s(d-a)}\over{4}}\right)$ are points of order 4 on $W$. Had we defined the birational equivalence between the projective closures of $W$ and $E$,  the points $\left(5d-a : \pm3s(d-a): 12\right)$ of $W$ would map to the singular point (0:1:0) of $E$, while the points $\left(-(a+d) \pm 6t: 0: 12\right)$ of $W$ would map to the singular point (1:0:0) of $E$.

\section{The Function Field of a Twisted Edwards Curve}

For Weierstrass curves 
$W: v^2=u^3+Au+B$
it is well known (see \cite{Silv} for example) that an
element of the function field $K(W)$ can be written uniquely in the form
\[
p(u)+vq(u)
\]
where $p(u) ,q(u)$ are polynomials in $u$.

We will prove an analogous result for twisted Edwards curves $E$.
Not surprisingly, rational functions are needed in place of the polynomials.
We use the notation $\text{ord}_P (f)$ to denote the valuation of a function
$f\in K(E)$ at a point $P$.

\begin{thm} \label{uniq1}
Any function $g \in K(E)$ can be written uniquely as \[ g(x, y) = p(y) + xq(y)\] where 
$p(y)$, $q(y)$ are rational functions in $y$.
\end{thm}

\emph{Proof:} 
Let $f(x,y)=0$ be the equation defining $E$, where \[ f(x,y) = ax^2+y^2-1-dx^2y^2. \]
In $K(E)$ we have
\[
x^2=\frac{1-y^2}{a-dy^2}.
\]
If $g(x,y)\in K(E)$, by replacing every occurence of $x^2$ by this rational function in $y$
it follows that $g(x,y)$ can be written in the form
\[
\frac{A(y)+xB(y)}{C(y)+xD(y)}
\]
where $A,B,C,D$ are rational functions.
Multiplying above and below by $C(y)-xD(y)$, and replacing each $x^2$ by
$\frac{1-y^2}{a-dy^2}$ shows that $g$ can be written in the stated form.
This proves existence.

Suppose for the sake of contradiction that this expression for $g$ is not unique. 
Then $A(y) + xB(y) = 0$ for some nonzero rational functions $A(y)$, $B(y)$. 
So \[ x = -{{A(y)}\over{B(y)}}\]  which implies
\begin{equation}\label{ords}
 \text{ord}_{(0,1)}x = \text{ord}_{(0,1)}A(y)-\text{ord}_{(0,1)}B(y).
 \end{equation}

 We obtain our contradiction by showing that the right-hand side of 
equation (\ref{ords}) is even,
but the left-hand side is equal to 1.

We expand at $(0,1)$ and we get 
\begin{align*}
f(x, y+1) &= ax^2+(y+1)^2-1-dx^2(y+1)^2\\
&=ax^2+y^2+2y-dx^2y^2-2dx^2y-dx^2.
\end{align*}
This shows that the line $x=0$ is not a tangent at $(0,1)$,
so $x$ is a local uniformizer there.
Then
\[ f(x, 0+1) = (a-d)x^2 \] 
which implies ord$_{(0,1)}(y-1) =2$ ord$_{(0,1)}(x) = 2$.


When computing ord$_{(0,1)}A(y)$, we translate $(0,1)$  to the origin, and write 
$A(y+1) = {{a(y)}\over{b(y)}}$ for some polynomials $a(y), \ b(y)$. 
Then 
\[\text{ord}_{(0,1)}A(y)= \text{ord}_{(0,0)}a(y)-\text{ord}_{(0,0)} b(y).\]
Of course, after translation we have $\text{ord}_{(0,0)}(y) = 2$.

Let 
$n_{0}$ be the degree of the term of smallest degree in $a(y)$,  
and similarly let $m_{0}$ be the degree of the term of smallest degree in $b(y)$. 
Then
$\text{ord}_{(0,0)}a(y) = \left(\text{ord}_{(0,0)}y\right) n_{0} = 2n_{0} $, and
similarly, $\text{ord}_{(0,0)}b(y) =  2m_{0}$. Thus $\text{ord}_{(0,1)}A(y)= 2(n_{0}-m_{0})$,
which is even.

Similarly,    $\text{ord}_{(0,1)}B(y)$ is even.
This proves that the right-hand side of (\ref{ords}) is even,
and we are done.
\hfill $\Box$
\bigskip

\begin{cor} \label{uniq2}
Any function $g \in K(E)$ can be written uniquely as 
\[ g(x, y) = p^{\prime}(y) + {{1}\over{x}}q^{\prime}(y)\] where 
$p^{\prime}(y)$, $q^{\prime}(y)$ are rational functions in $y$.
\end{cor}

\emph{Proof:} 
This follows from the Theorem \ref{uniq1}, and
the fact that \[x={{1}\over{x}}\cdot {{1-y^2}\over{a-dy^2}}\] on the function field of $E$.
In fact $p^\prime (y)$ is equal to $p(y)$, using the notation of Theorem \ref{uniq1},
and 
\[
q^\prime (y) = {{1-y^2}\over{a-dy^2}} \ q(y).
\]
\hfill $\Box$

\section{Division Rational Functions on Twisted Edwards Curves}

We define the following rational functions $\psi_n(x,y)$ on the function field of $E$ 
recursively for $n\geq0$:

\begin{align*}\psi_0(x,y) &:= 0\\
\psi_1(x,y) &:= 1\\
\psi_2(x,y) &:= {{(a-d)(1+y)}\over{x(2(1-y))}}\\
\psi_3(x,y) &:= {{(a-d)^3(a+2ay-2dy^3-dy^4)}\over{(2(1-y))^4}}\\
\psi_4(x,y) &:= {{2(a-d)^6y(1+y)(a-dy^4)}\over{x((2(1-y))^7}}\\
\psi_{2m+1}(x,y) &:= \psi_{m+2}(x,y)\psi_m^3(x,y) - \psi_{m-1}(x,y)\psi_{m+1}^3(x,y) \quad \text{for } m \geq 2\\
\psi_{2m}(x,y) &:= {\psi_m(x,y)\over\psi_2(x,y)}\left(\psi_{m+2}(x,y)\psi_{m-1}^2(x,y) - \psi_{m-2}(x,y)\psi_{m+1}^2(x,y)\right) \quad \text{for } m \geq 3.\end{align*}

These functions are not defined at the points $(0,1)$ and $(0,-1)$.
We point out that these elements of the function field $K(E)$ are 
in the unique form given in Corollary \ref{uniq2}.

For $n \geq 1$, we also define

\begin{align*}\phi_n(x,y) &: = {{(1+y)\psi_n^2(x,y)}\over{(1-y)}} - {{4\psi_{n-1}(x,y)\psi_{n+1}(x,y)}\over{(a-d)}}\\
\text{and} \quad\omega_n(x,y) &:= {{2\psi_{2n}(x,y)}\over{(a-d)\psi_n(x,y)}}.\end{align*}

Next we show that these rational functions arise in the
multiplication-by-$n$ map.

\begin{thm}\label{divform}
Let $(x,y)$ be a point in $E(\overline{k})\setminus \{ (0,1), (0,-1)\} $ and $n \ge 1$ an integer. Then \[[n](x,y) = \left({{\phi_n(x,y)\psi_n(x,y)}\over{\omega_n(x,y)}}, {{\phi_n(x,y) - \psi_n^2(x,y)}\over{\phi_n(x,y) + \psi_n^2(x,y)}}\right).\]
\end{thm}

\emph{Proof:} 
Compute the division polynomials for the  Weierstrass  elliptic curve from Section \ref{tecsection},
$W: v^2 = u^3 + Au + B, \ $ where

\[A = - {{(a^2+14ad+d^2)}\over{48}},\qquad B = - {{(a^3-33a^2d-33ad^2+d^3)}\over{864}}.\]
We get
\begin{align*}\Psi_{0}(u, v) &= 0\\
\Psi_{1}(u, v) &= 1\\
\Psi_{2}(u, v) &= 2v\\
\Psi_{3}(u, v) &= 3u^4+6Au^2+12Bu - A^2\\
\Psi_{4}(u, v) &= 4v(u^6+5Au^4+20Bu^3-5A^2u^2-4ABu-A^3-8B^2)\\
\Psi_{2m+1}(u, v) &= \Psi_{m+2}(u,v)\Psi_{m}^3(u,v) - \Psi_{m-1}(u,v)\Psi_{m+1}^3(u,v) \text{ for } m \geq 2\\
\Psi_{2m}(u,v) &= {\Psi_{m}(u,v)\over\Psi_{2}(u,v)}\left(\Psi_{m+2}(u,v)\Psi_{m-1}^2(u,v) - \Psi_{m-2}(u,v)\Psi_{m+1}^2(u,v)\right)\quad \text{ for }m \geq 3.\end{align*}

Substituting \begin{align*}&A = - {{(a^2+14ad+d^2)}\over{48}},\quad B = - {{(a^3-33a^2d-33ad^2+d^3)}\over{864}}\quad \text{and}\\
&u := {{(5a-d)+(a-5d)y}\over{12(1-y)}},   \qquad v := {{(a-d)(1+y)}\over{4x(1-y)}},\end{align*}

for the cases $0,1,2,3,4$ we see that $\Psi_{i}(u, v) = \psi_i(x,y)$ for $i = 0,1,2,3,4$. Hence, as the recursion relations for the two sets of functions $\Psi_{i}(u, v)$ and $\psi_i(x,y)$ 
are identical for $i \geq 5$, we have that $\Psi_{n}(u, v) = \psi_n(x,y)$ for all integers $n \geq 0$.

From here on we will use the abbreviated notations $\psi_n$ for $\psi_n(x,y)$, $\phi_n$ for $\phi_n(x,y)$ and $\omega	_n$ for $\omega_n(x,y)$.
Let $\left(x_n , y_n\right) = [n](x,y)$, and $\left(u_n, v_n\right) = [n]_{W}\left(u, v\right)$.

From the properties of the division polynomials,

\[u_n = u -{{\Psi_{n-1}(u, v)\Psi_{n+1}(u, v)}\over{\Psi_{n}^2(u, v)}}, \quad
v_n = {{\Psi_{2n}(u, v)}\over{2\Psi_{n}^4(u, v)}},\]

i.e.,

\[u_n = u -{{\psi_{n-1}\psi_{n+1}}\over{\psi_{n}^2}}, \quad
v_n = {{\psi_{2n}}\over{2\psi_{n}^4}},\]

and, applying the birational equivalence gives

\[x_n = {{6u_n - (a+d)}\over{6v_n}}, \quad y_n =  {{12u_n+d-5a}\over{12u_n+a-5d}},\]

\begin{align*}x_n &= {{2\psi_n^4}\over{\psi_{2n}}}\left({{5a-d+(a-5d)y}\over{12(1-y)}} - {{\psi_{n-1}\psi_{n+1}}\over{\psi_n^2}}- {{a+d}\over{6}}\right)\\
&= {{\psi_n^2}\over{\psi_{2n}}}\left({{(a-d)(1+y)\psi_n^2}\over{2(1-y)}}- {2\psi_{n-1}\psi_{n+1}}\right)\end{align*}

while \begin{align*} {{\phi_n\psi_n}\over{\omega_n}} &= {{(a-d)\psi_n^2}\over{2\psi_{2n}}}\left(\left({{1+y}\over{1-y}}\right)\psi_n^2-{{4\psi_{n-1}\psi_{n+1}}\over{a-d}}\right)\\
&= {{\psi_n^2}\over{\psi_{2n}}}\left({{(a-d)(1+y)\psi_n^2}\over{2(1-y)}}- {2\psi_{n-1}\psi_{n+1}}\right)\\
 &=x_n.\\
\end{align*}
Also,
\[y_n = {{12u_n+d-5a}\over{12u_n+a-5d}}\]
and
\begin{align*}12u_n+d-5a &= {{5a-d+(a-5d)y}\over{(1-y)}} - 12{{\psi_{n-1}\psi_{n+1}}\over{\psi_n^2}} +d-5a\\
&={{6(a-d)y}\over{1-y}}- 12{{\psi_{n-1}\psi_{n+1}}\over{\psi_n^2}}\\
12u_n+a-5d &= {{6(a-d)}\over{1-y}}- 12{{\psi_{n-1}\psi_{n+1}}\over{\psi_n^2}}\end{align*}

so \[y_n = {{(a-d)y\psi_n^2-2(1-y)\psi_{n-1}\psi_{n+1}}\over{(a-d)\psi_n^2-2(1-y)\psi_{n-1}\psi_{n+1}}}\]

and \begin{align*}{{\phi_n-\psi_n^2}\over{\phi_n+\psi_n^2}} &= {{\left({{1+y}\over{1-y}}\right)\psi_n^2-{{4\psi_{n-1}\psi_{n+1}}\over{a-d}}-\psi_n^2}\over{\left({{1+y}\over{1-y}}\right)\psi_n^2-{{4\psi_{n-1}\psi_{n+1}}\over{a-d}}+\psi_n^2}}\\
&={{(a-d)y\psi_n^2-2(1-y)\psi_{n-1}\psi_{n+1}}\over{(a-d)\psi_n^2-2(1-y)\psi_{n-1}\psi_{n+1}}}\\
& = y_n .\end{align*}

Hence \[[n](x,y) = \left({{\phi_n(x,y)\psi_n(x,y)}\over{\omega_n(x,y)}}, {{\phi_n(x,y) - \psi_n^2(x,y)}\over{\phi_n(x,y) + \psi_n^2(x,y)}}\right).\] 
\hfill $\Box$

\begin{cor}\label{divformcor}
Let $P=(x,y)$ be in
$E(\overline{k})\setminus \{ (0,1), (0,-1)\} $ and let $n\geq 1$.
Then $P$ is an $n$-torsion point of $E$ if and 
only if $\psi_n(P) = 0$.
\end{cor}
  
\emph{Proof:}  Since the identity is $(0,1)$,
the result is clear from Theorem \ref{divform}.
\hfill $\Box$

\bigskip

 So the $\psi_n(x,y)$, though they are rational functions, can be seen as analogues of division polynomials.  Here are the first seven $\psi_n(x,y)$:

\begin{align*}
\psi_0 &= 0\\
\psi_1 &= 1\\
\psi_2 &= {{(a-d)(y+1)}\over{x(2(1-y))}}\\
\psi_3 &= {{(a-d)^3(-dy^4-2dy^3+2ay+a)}\over{(2(1-y))^4}}\\
\psi_4 &= {{2(a-d)^6(-dy^6-dy^5+ay^2+ay)}\over{x((2(1-y))^7}}\\
\psi_5 &={{(a-d)^9(d^3y^{12}-2d^3y^{11}+\dots+2a^3y-a^3)}\over{(2(1-y))^{12}}}\\
\psi_6 &= {{(a-d)^{13}(-d^4y^{17}-d^4y^{16}+(4ad^3+4d^4)y^{15}+\dots+(4a^3d+4a^4)y^2-a^4y-a^4)}\over{x((2(1-y))^{17}}}.
\end{align*}

As we said earlier,  these elements of the function field $K(E)$ are 
in the unique form given in Corollary \ref{uniq2}.

The apparent patterns here are proved in theorem \ref{main} below.

\section{Gauss's notes} \label{gauss}

We mention here how Gauss's formulas (see Fig 1) are incorrect,
although they are close to being correct.
Essentially the only errors are sign errors.


One can see that Gauss calls the point  $(s,c)$ and sin lemn $n\varphi$ denotes the
$x$ coordinate of   $[n](s,c)$, and cos lemn $n\varphi$ denotes the
$y$ coordinate of   $[n](s,c)$.

 We represent our formulas in the unique form given by Theorem \ref{uniq1}.
 
Our division polynomial formulas applied to the curve (\ref{gauss1}) give
\[
[2](s,c) = \biggl(\frac{2sc ( c^2 +  1)}{c^4 + 1}, \frac{-c^4 - 2 c^2 + 1}{c^4 - 2 c^2 - 1}\biggr)
\]
which we can see agree with Gauss's formula for twice the point in terms of $c$. However, there is an error in Gauss's formula for cos lemn $2 \varphi$ in terms of $s$, which should be \[\frac{1-2s^2-s^4}{1+2s^2-s^4}\]

A sign error also occurs in the denominator of the sin lemn $5 \varphi$ formula (coefficient of $s^{12}$ should be -12), and six times in the cos lemn $4 \varphi$ formula,
which should read
\[\text{cos lemn}4 \varphi = \frac{1-8s^2-12s^4-8s^6++38s^8+8s^{10}-12s^{12}+8s^{14}+s^{16}}{1+8s^2-12s^4+8s^6+38s^8-8s^{10}-12s^{12}-8s^{14}+s^{16}}\]

We note that these sign errors break the apparent ``reverse symmetry'' between the coefficients of the  numerator and denominator. This symmetry, proved by Abel \cite{Abel}, is explained in greater detail in Chapter 15 of \cite{Cox}. 

For the general case, Gauss gave some information on
the $x$ coordinate of $[n](s,c)$, but not the $y$ coordinate.

\section{Division Polynomials} \label{divpols}

The next theorem isolates the key polynomial in the numerator of $\psi_n$,
which we call $\tilde{\psi} (y)$.
These polynomials could also be called the division polynomials for twisted Edwards curves.

\begin{thm} \label{main}
We have
\begin{equation*}
\psi_n(x,y)= \left\{
	\begin{array}{rl}
		{(a-d)^{k(n)}\tilde{\psi}_{n}(y)/(2(1-y))^{m(n)}} & \text{ if $n$ is odd}\\
		\\
		{(a-d)^{k(n)}\tilde{\psi}_{n}(y)/x(2(1-y))^{m(n)}} & \text{ if $n$ is even}
	\end{array} \right.
\end{equation*}
where
\begin{equation*}
m(n)= \left\{
	\begin{array}{rl}
		{{n^2-1}\over{2}} & \text{ if n is odd}\\
		\\
		{{n^2-2}\over{2}} & \text{ if n is even}
	\end{array} \right.
\end{equation*}
and
\[
k(n) = \left\lfloor{{3n^2}\over{8}}\right\rfloor\]
and
\begin{align*}
\tilde{\psi}_0(y) &= 0\\
\tilde{\psi}_1(y) &= 1\\
\tilde{\psi}_2(y) &= y+1\\
\tilde{\psi}_3(y) &= -dy^4-2dy^3+2ay+a\\
\tilde{\psi}_4(y) &= -2y(y+1)(dy^4-a) = -2dy^6-2dy^5+2ay^2+2ay,\end{align*}
and
\begin{equation*}
\tilde{\psi}_{2r+1}(y) = \left\{
	\begin{array}{rl}
		{{4(a-d)(a-dy^2)^2\tilde{\psi}_{r+2}(y)\tilde{\psi}_{r}^3(y)}\over{(y+1)^2}} - \tilde{\psi}_{r-1}(y)\tilde{\psi}_{r+1}^3(y)& \text{ if } r\equiv 0 \pmod{4},\ r\geq4\\
		{\tilde{\psi}_{r+2}(y)\tilde{\psi}_{r}^3(y)} - {{4(a-dy^2)^2\tilde{\psi}_{r-1}(y)\tilde{\psi}_{r+1}^3(y)}\over{(y+1)^2}}& \text{ if } r\equiv 1 \pmod{4},\ r\geq5\\
		{{4(a-dy^2)^2\tilde{\psi_{r+2}}(y)\tilde{\psi_{r}}^3(y)}\over{(y+1)^2}} - \tilde{\psi}_{r-1}(y)\tilde{\psi}_{r+1}^3(y)& \text{ if } r\equiv 2 \pmod{4},\ r\geq2\\
		{\tilde{\psi}_{r+2}(y)\tilde{\psi}_{r}^3(y)} - {{4(a-d)(a-dy^2)^2\tilde{\psi}_{r-1}(y)\tilde{\psi}_{r+1}^3(y)}\over{(y+1)^2}}& \text{ if } r\equiv 3 \pmod{4},\ r\geq3
	\end{array} \right.
\end{equation*}
and
\begin{equation*}
\tilde{\psi}_{2r}(y) = \left\{
	\begin{array}{rl}
		{{\tilde{\psi}_{r}(y)}\over{y+1}}\left(\tilde{\psi}_{r+2}(y)\tilde{\psi}_{r-1}^2(y) - \tilde{\psi}_{r-2}(y)\tilde{\psi}_{r+1}^2(y)\right)& \text{ if } r\equiv 0 \pmod{4},\ r\geq4\\
		{{\tilde{\psi}_{r}(y)}\over{y+1}}\left((a-d)\tilde{\psi}_{r+2}(y)\tilde{\psi}_{r-1}^2(y) - \tilde{\psi}_{r-2}(y)\tilde{\psi}_{r+1}^2(y)\right)& \text{ if } r\equiv 1 \pmod{4},\ r\geq5\\
		{{\tilde{\psi}_{r}(y)}\over{y+1}}\left(\tilde{\psi}_{r+2}(y)\tilde{\psi}_{r-1}^2(y) - \tilde{\psi}_{r-2}(y)\tilde{\psi}_{r+1}^2(y)\right)& \text{ if } r\equiv 2 \pmod{4},\ r\geq6\\
		{{\tilde{\psi}_{r}(y)}\over{y+1}}\left(\tilde{\psi}_{r+2}(y)\tilde{\psi}_{r-1}^2(y) - (a-d)\tilde{\psi}_{r-2}(y)\tilde{\psi}_{r+1}^2(y)\right)& \text{ if } r\equiv 3 \pmod{4},\ r\geq3.
		\end{array} \right.
\end{equation*}
\end{thm}

\emph{Proof:} 

First observe for all $ t \in \mathbb{Z}, \ t>0$,
\begin{align*}m(4t) &= {{16t^2-2}\over{2}} = 8t^2-1\\
m(4t\pm1)&={{(4t\pm1)^2-1}\over{2}}={{16t^2\pm8t}\over{2}}=8t^2\pm4t\\
m(4t\pm2)&={{(4t\pm2)^2-2}\over{2}}={{16t^2\pm16t+2}\over{2}} = 8t^2\pm8t+1\\
m(4t\pm3)&={{(4t\pm3)^2-1}\over{2}}={{16t^2\pm24t+8}\over{2}} = 8t^2\pm12t+4
\end{align*}
and
\begin{align*}
k(4t) &= \left\lfloor{{3(4t)^2}\over{8}}\right\rfloor = \left\lfloor{6t^2}\right\rfloor = 6t^2\\
k(4t\pm1) &= \left\lfloor{{{3(4t\pm1)^2}\over{8}}}\right\rfloor =\left\lfloor{6t^2\pm3t+{{3}\over{8}}}\right\rfloor= 6t^2\pm3t\\
k(4t\pm2) &= \left\lfloor{{{3(4t\pm2)^2}\over{8}}}\right\rfloor =\left\lfloor{6t^2\pm6t+{{12}\over{8}}}\right\rfloor= 6t^2\pm6t+1\\
k(4t\pm3) &= \left\lfloor{{{3(4t\pm3)^2}\over{8}}}\right\rfloor =\left\lfloor{6t^2\pm9t+{{27}\over{8}}}\right\rfloor= 6t^2\pm9t+3.
\end{align*}

The proof is by induction.
The claim is true for $n = 0 \dots 4$.

Assume true for $0 \dots n-1$

\underline{Case 1:} $n \equiv 0 \pmod{8}$ i.e. $ n = 8l$ for some $l \in \mathbb{Z}$. Let $r = 4l$.

By definition,
\begin{align*}\psi_{n} &= {{\psi_r}\over{\psi_2}}\left(\psi_{r+2}\psi_{r-1}^2 - \psi_{r-2}\psi_{r+1}^2\right)\\
&= {{(a-d)^{k(r)-1}\tilde{\psi}_{r}}\over{(y+1)(2(1-y))^{m(r)-1}}}\left({{(a-d)^{k(r+2)+2k(r-1)}\tilde{\psi}_{r+2}\tilde{\psi}_{r-1}^2}\over{x(2(1-y))^{m(r+2)+2m(r-1)}}} - {{(a-d)^{k(r-2)+2k(r+1)}\tilde{\psi}_{r-2}\tilde{\psi}_{r+1}^2}\over{x(2(1-y))^{m(r-2)+2m(r+1)}}}\right).\end{align*}

Also,
\begin{align*}m(4l)-1+m(4l+2) + 2m(4l-1) &= 8l^2-1-1+8l^2+8l+1+16l^2-8l\\
&=32l^2-1 = m(8l) = m(n)\\
m(4l)-1+m(4l-2) + 2m(4l+1) &= 8l^2-1-1+8l^2-8l+1+16l^2+8l\\
&=32l^2-1 = m(8l) = m(n)\end{align*}

and
\begin{align*}k(4l)-1+k(4l+2) + 2k(4l-1) &= 6l^2-1+6l^2+6l+1+12l^2-6l\\
&=24l^2 = k(8l) = k(n)\\
k(4l)-1+k(4l-2) + 2k(4l+1) &= 6l^2-1+6l^2-6l+1+12l^2+6l\\
&=24l^2 = k(8l) = k(n).
\end{align*}

So \begin{align*}
\psi_n &={{(a-d)^{k(n)}}\over{x(y+1)(2(1-y))^{m(n)}}}\left(\tilde{\psi}_{r}\left(\tilde{\psi}_{r+2}\tilde{\psi}_{r-1}^2 - \tilde{\psi}_{r-2}\tilde{\psi}_{r+1}^2\right)\right)\\
&={{(a-d)^{k(n)}\tilde{\psi}_{n}(y)}\over{x(2(1-y))^{m(n)}}} \ .
\end{align*}

\underline{Case 2:}  $n \equiv 1 \pmod{8}$ i.e. $ n = 8l+1$ for some $l \in \mathbb{Z}$. Let $r = 4l$.

By definition
\begin{align*}\psi_{n} &= \psi_{r+2}\psi_r^3 - \psi_{r-1}\psi_{r+1}^3\\
&={{(a-d)^{k(r+2)+3k(r)}\tilde{\psi}_{r+2}\tilde{\psi}_r^3}\over{y^4(2(1-x))^{m(r+2)+3m(r)}}} - {{(a-d)^{k(r-1)+3k(r+1)}\tilde{\psi}_{r-1}\tilde{\psi}_{r+1}^3}\over{(2(1-y))^{m(r-1)+3m(r+1)}}}.
\end{align*}

Using the curve equation \[ax^2+y^2  = 1 + dx^2y^2\]
gives\begin{align*}&x^2 = {{(1-y^2)}\over{(a-dy^2)}} = {{(1-y)(1+y)}\over{(a-dy^2)}}\\
\Rightarrow &x^4= {{(1-y)^2(1+y)^2}\over{(a-dy^2)^2}}\end{align*}

so \[\psi_n = {{4(a-d)^{k(r+2)+3k(r)}(a-dy^2)^2\tilde{\psi}_{r+2}\tilde{\psi}_r^3}\over{(y+1)^2(2(1-y))^{m(r+2)+3m(r)+2}}} - {{(a-d)^{k(r-1)+3k(r+1)}\tilde{\psi}_{r-1}\tilde{\psi}_{r+1}^3}\over{(2(1-y))^{m(r-1)+3m(r+1)}}}\ .\]

Again,
\begin{align*}k(4l+2)+3k(4l) & =6l^2+6l+1 + 18l^2 = 24l^2+6l+1\\& = k(n)+1\\
k(4l-1)+3k(4l+1)&= 6l^2-3l+18l^2+9l=24l^2+6l \\&= k(n)
\end{align*}

and
\begin{align*}
m(4l+2)+3m(4l)+2 &= 8l^2+8l+1+24l^2-3 +2= 32l^2 +8l\\&=m(n)\\
m(4l-1)+3m(4l+1) &=8l^2-4l+24l^2+12l = 32l^2 +8l\\&=m(n).
\end{align*}
Hence\[\psi_n = {{4(a-d)(a-dy^2)^2\tilde{\psi}_{r+2}(y)\tilde{\psi}_{r}^3(y)}\over{(y+1)^2}} - \tilde{\psi}_{r-1}(y)\tilde{\psi}_{r+1}^3(y)\ .\]

\underline{Cases 3,$\dots$8:}  $n \equiv 2, \dots 7 \pmod{8}$. Similar.
\hfill $\Box$

\bigskip

\begin{cor}\label{coronevar}  
Let $P=(x,y)$ be in
$E(\overline{k})\setminus \{ (0,1)\} $ and let $n\geq 1$.
Then 
\begin{center}
$P$ is an $n$-torsion point of $E$\quad  if and only if  \quad  $\tilde\psi_n(y) = 0$.
\end{center}
\end{cor}
  
\emph{Proof:} 
The result follows from Corollary \ref{divformcor} and Theorem \ref{main}.
\hfill $\Box$

\bigskip

\section{Further Facts}\label{further}

Here are some more facts about the $\tilde{\psi}$.

\begin{thm} \label{int}
$\tilde{\psi}_n(y) \in \mathbb{Z}[a,d,y] \ \ \forall n >0$, and $(y+1)$ divides $\tilde{\psi}_n(y)$ if n is even
\end{thm}

\emph{Proof:} Proof is by induction. The statement is true for $n = 0, 1, 2, 3, 4$. Now suppose it is true for $0,1,2, \dots, n-1$:

\underline{Case 1:} $n \equiv 0 \pmod{8}$ i.e. $ n = 8l$ for some $l \in \mathbb{Z}$. Let $r = 4l$.

Then $\tilde{\psi}_n(y) = {{\tilde{\psi}_{r}(y)}\over{y+1}}\left(\tilde{\psi}_{r+2}(y)\tilde{\psi}_{r-1}^2(y) - \tilde{\psi}_{r-2}(y)\tilde{\psi}_{r+1}^2(y)\right)$

and $\tilde{\psi}_{r}(y), \ \tilde{\psi}_{r+2}(y), \ \tilde{\psi}_{r-1}(y),\ \tilde{\psi}_{r-2}(y), \ \tilde{\psi}_{r+1}(y) \in \mathbb{Z}[a,d,y] $. Also, $(y+1)$ divides $\tilde{\psi}_{r}(y), \ \tilde{\psi}_{r+2}(y)$, and $\tilde{\psi}_{r-2}(y)$ by hypothesis. Hence $\tilde{\psi}_n(y) \in \mathbb{Z}[a,d,y] $ and $(y+1)$ divides $\tilde{\psi}_{n}(y)$.

\underline{Case 2:}  $n \equiv 1 \pmod{8}$ i.e. $ n = 8l+1$ for some $l \in \mathbb{Z}$. Let $r = 4l$.

Then $\tilde{\psi}_n(y) = {{4(a-d)(a-dy^2)^2\tilde{\psi}_{r+2}(y)\tilde{\psi}_{r}^3(y)}\over{(y+1)^2}} - \tilde{\psi}_{r-1}(y)\tilde{\psi}_{r+1}^3(y)$

and $\tilde{\psi}_{r+2}(y), \ \tilde{\psi}_{r}(y), \ \tilde{\psi}_{r-1}(y), \ \tilde{\psi}_{r+1}(y) \in \mathbb{Z}[a,d,y] $. Also, $(y+1)$ divides $\tilde{\psi}_{r}(y)$ and $\tilde{\psi}_{r+2}(y)$ by hypothesis. Hence $\tilde{\psi}_n(y) \in \mathbb{Z}[a,d,y]$.

\underline{Cases 3,$\dots$8:}  $n \equiv 2, \dots 7 \pmod{8}$. Similar.
\hfill $\Box$
\bigskip

Theorem \ref{leading} and Corollary \ref{degree} provide results for the degrees of these polynomials $\tilde{\psi}_n(y)$, and  Theorem \ref{symm} shows that the coefficients of the polynomials exhibit a large amount of symmetry.

\begin{thm} \label{leading} If $char(k) = 0$ or $4 \cdot char(k) \nmid n$, then $\tilde{\psi}_{n}(y)$ has leading term (term of largest degree in $y$)
\begin{equation*} \left\{
	\begin{array}{rl}
		\delta(n) d^{m(n)-k(n)}y^{m(n)}& \text{ if } n \not\equiv 0 \pmod{4}\\
		\\
		\delta(n) d^{m(n)-k(n)}y^{m(n)-1}& \text{ if } n \equiv 0 \pmod{4}
	\end{array} \right.
\end{equation*}
where
\begin{equation*} \delta (n) = \left\{
	\begin{array}{rl}
		{{n}\over{2}} & \text{ if } n \equiv 0 \pmod{8}\\
		\\
		- {{n}\over{2}} & \text{ if } n \equiv 4 \pmod{8}\\
		\\
		1 & \text{ if } n \equiv 1,2, \text{ or } 5 \pmod{8}\\
		\\
		-1& \text{ if } n \equiv 3,6, \text{ or } 7 \pmod{8}
	\end{array} \right.
\end{equation*}
and $m(n)$, $k(n)$ are as defined in Theorem \ref{main}.

If $char(k) \neq 0$ and $4\cdot char(k) \mid n$, then $deg(\tilde{\psi}_{n}(y)) < m(n)-1$ .
\end{thm}

\emph{Proof:} Proof is by induction. The statement is true for  $n = 0, 1, 2, 3, 4$. Now suppose it is true for $0,1,2, \dots, n-1$:

\underline{Case 1:} $n \equiv 0 \pmod{8}$ i.e. $ n = 8l$ for some $l \in \mathbb{Z}$. Let $r = 4l$.
Then
\begin{align*}
	\tilde{\psi}_n(y) =& {{\tilde{\psi}_{r}(y)}\over{y+1}}\left(\tilde{\psi}_{r+2}(y)\tilde{\psi}_{r-1}^2(y) - \tilde{\psi}_{r-2}(y)\tilde{\psi}_{r+1}^2(y)\right)\\
	=&(\delta(r) d^{m(r)-k(r)}y^{m(r)-2}+\dots)\times \\
	&[(\delta(r+2)(\delta(r-1))^2d^{m(r+2)+2m(r-1)-k(r+2)-2k(r-1)}y^{m(r+2)+2m(r-1)} +\dots)\\& \  -(\delta(r-2)(\delta(r+1))^2d^{m(r-2)+2m(r+1)-k(r-2)-2k(r+1)}y^{m(r-2)+2m(r+1)}+ \dots)]
\end{align*}

So, computing the $m$'s and $k$'s as in previous proofs, and noting that \begin{align*}\delta(r) = \pm2l,\ &\delta(r+2) = \pm1, \ \delta(r-1) = -1, \\ &\delta(r-2) = \mp1, \ \delta(r+1) = 1, \end{align*}
the leading term is thus 
\[ \pm2ld^{m(n)-k(n)}y^{m(r)-2}(\pm y^{m(r+2)+2m(r-1)} \pm y^{m(r-2)+2m(r+1)})\] \[ = {{n}\over{2}}d^{m(n)-k(n)}y^{m(n)-1} \]
\[=\delta(n)d^{m(n)-k(n)}y^{m(n)-1}. \]

The only exception being if $char(k) \neq 0$ and $char(k) \mid r$, (i.e. if $char(k) \mid n$) in which case,  $deg(\tilde{\psi}_{r}(y)) < m(r)-1$ and  $deg(\tilde{\psi}_{n}(y)) < m(n)-1$.

\underline{Case 2:}  $n \equiv 1 \pmod{8}$ i.e. $ n = 8l+1$ for some $l \in \mathbb{Z}$. Let $r = 4l$.

Then $\tilde{\psi}_n(y) = {{4(a-d)(a-dy^2)^2\tilde{\psi}_{r+2}(y)\tilde{\psi}_{r}^3(y)}\over{(y+1)^2}} - \tilde{\psi}_{r-1}(y)\tilde{\psi}_{r+1}^3(y)$.

The degree (in $y$) of the first term above is $m(r+2)+3(m(r)-1)+4-2 = 32l^2+8l-3$.

The degree (in $y$) of the second term is $m(r-1)+3m(r+1) = 32l^2+8l$
Thus ${{4(a-d)(a-dy^2)^2\tilde{\psi}_{r+2}(y)\tilde{\psi}_{r}^3(y)}\over{(y+1)^2}}$ does not contribute to the leading term which is 
\[-\delta(r-1)(\delta(r+1))^3d^{m(r-1)+3m(r+1)-k(r-1)-3k(r+1)}y^{32l^2+8l}. \] Now,
\[\delta(r-1) = -1, \  \delta(r+1) = 1, \ \delta(n) = 1\]
\[k(r-1)+3k(r+1) =24l^2+6l\]
\[m(n) = m(8l+1)=32l^2+8l-(24l^2+6l) = 8l^2+2l. \]
So the leading term is $d^{m(n)-k(n)}y^{m(n)} = \delta(n)d^{m(n)-k(n)}y^{m(n)}$, as required.

The only exceptional case is if $char(k) \neq 0$ and $char(k) \mid r$, in which case $deg(\tilde{\psi}_{r}(y)) < m(r)-1$, but as  $\tilde{\psi}_{r}(y)$ does not contribute to the leading term, this does not affect the result.

\underline{Cases 3,$\dots$8:}  $n \equiv 2, \dots 7 \pmod{8}$. Similar.
\hfill $\Box$
\bigskip

\begin{cor} \label{degree}  If $4 \nmid n$, then $deg(\tilde{\psi}_{n}(y)) = m(n)$ where
\begin{equation*}
m(n)= \left\{
	\begin{array}{rl}
		{{n^2-1}\over{2}} & \text{ if n is odd}\\
		\\
		{{n^2-2}\over{2}} & \text{ if n is even.}
	\end{array} \right.
\end{equation*}
If $4 \mid n$ and $char(k) \nmid n$,  $deg(\tilde{\psi}_{n}(y)) = m(n)-1.$

Otherwise $deg(\tilde{\psi}_{n}(y))< m(n)-1.$
\end{cor}

\emph{Proof:} Immediate from Theorem \ref{leading} .   \hfill $\ \Box$

The only case where the degree of the polynomial $\tilde{\psi}_{n}$ is not known precisely is when $4 \cdot char(k) \mid n$. In any case, $\frac{n^2}{2}$ is an upper bound for $deg(\tilde{\psi}_{n})$.

\begin{lem} \label{last} If $char(k) = 0$ or $4 \cdot char(k) \nmid n$, then $\tilde{\psi}_{n}(y)$ has final term (term of least degree in $y$)
\begin{equation*} \left\{
	\begin{array}{rl}
		\epsilon(n) a^{m(n)-k(n)}& \text{ if } n \not\equiv 0 \pmod{4}\\
		\\
		\epsilon(n) a^{m(n)-k(n)}y& \text{ if } n \equiv 0 \pmod{4}
	\end{array} \right.
\end{equation*}
where
\begin{equation*} \epsilon (n) = \left\{
	\begin{array}{rl}
		-{{n}\over{2}} & \text{ if } n \equiv 0 \pmod{8}\\
		\\
		 {{n}\over{2}} & \text{ if } n \equiv 4 \pmod{8}\\
		\\
		1 & \text{ if } n \equiv 1,2, \text{ or } 3 \pmod{8}\\
		\\
		-1& \text{ if } n \equiv 5,6, \text{ or } 7 \pmod{8}
	\end{array} \right.
\end{equation*}
and $m(n)$, $k(n)$ are as defined in Theorem \ref{main}.

If $char(k) \neq 0$ and $4\cdot char(k) \mid n$, then the term of least degree has degree greater than 1. 
\end{lem}

\emph{Proof:} Similar to proof of Theorem \ref{leading}.  
\hfill $\ \Box$
\bigskip

Recall from Theorem \ref{int} that $\tilde{\psi}_{n}(y) = \tilde{\psi}_{n}(a,d,y) \in  \mathbb{Z}[a,d,y]$. If we write $\tilde{\psi}_{n}$ in the form \[ \tilde{\psi}_{n}(a,d,y) = \alpha_{m(n)}y^{m(n)} + \alpha_{m(n)-1}y^{m(n)-1}+ \dots + \alpha_{1}y+\alpha_{0}\] where $m(n)$ is as defined in Theorem \ref{main} (so, in particular, if $4 \mid n, \  \alpha_{m(n)} = \alpha_{0} =0$) and $\alpha_{i} \in \mathbb{Z}[a,d]$, then we define \[\tilde{\psi}_{n}^*(a,d,y):=\alpha_{0}y^{m(n)}+\alpha_{1}y^{m(n)-1}+\dots+\alpha_{m(n)-1}y+\alpha_{m(n)}\]

\begin{lem} \label{hom} $\tilde{\psi}_{n}(a,d,y)$, considered as a polynomial in $a$ and $d$ (with 
coefficients in $\mathbb{Z}[a,d]$) is homogeneous of degree $m(n)-k(n)$.
\end{lem}

\emph{Proof:} Proof is by induction using Theorem \ref{main}.  
\hfill $\ \Box$
\bigskip

\begin{thm} \label{symm}
Consider $\tilde{\psi}_{n}(a,d,y) \in  \mathbb{Z}[a,d,y]$, as a polynomial in three variables. Then $\tilde{\psi}_{n}(a,d,y) = \tilde{\psi}_{n}^*(-d,-a,y)$.
\end{thm}

\emph{Proof:} We can restate this theorem as: If \[\tilde{\psi}_{n}(a,d,y) = \alpha_{m(n)}(a,d)y^{m(n)} + \alpha_{m(n)-1}(a,d)y^{m(n)-1}+ \dots + \alpha_{1}(a,d)y+\alpha_{0}(a,d)\]then \[\tilde{\psi}_{n}(a,d,y) = \alpha_{0}(-d,-a)y^{m(n)}+\alpha_{1}(-d,-a)y^{m(n)-1}+\dots+\alpha_{m(n)-1}(-d,-a)y+\alpha_{m(n)}(-d,-a).\]

If $E$ is as defined at the outset, \[E: ax^2+y^2  = 1+dx^2y^2\]
and we let $E^{\prime}$ be the twisted Edwards curve \[E^{\prime}: dx^2+y^2 = 1+ax^2y^2\]
then the birational equivalence $(x,y) \mapsto \left(x, {{1}\over{y}}\right)$ maps $E$ to $E^{\prime}$, and $E^{\prime}$ to $E$.

Now, \[\psi_{n}(x,y) = {{(a-d)^{k(n)}\tilde{\psi}_{n}(y)}\over{(2(1-y))^{m(n)}x^{\gamma(n)}}}\] where
 \begin{equation*} \gamma(n) = \left\{
	\begin{array}{rl}
		1 &\text{ if $n$ is even}\\
		0 &\text{ if $n$ is odd}
	\end{array} \right.
\end{equation*}

and \[\psi_{n}^{\prime}(x,y) = {{(d-a)^{k(n)}\tilde{\psi}_{n}^{\prime}(y)}\over{(2(1-y))^{m(n)}x^{\gamma(n)}}}\]

where $\psi_{n}^{\prime}(x,y), \ \tilde{\psi}_{n}^{\prime}(y)$ are the relevant functions defined on $E^{\prime}$.

Now,
\begin{align*}
\psi_{n}^{\prime}(x, {{1}\over{y}}) &= {{(d-a)^{k(n)}\tilde{\psi}_{n}^{\prime}({{1}\over{y}})}\over{(2(1-{{1}\over{y}}))^{m(n)}x^{\gamma(n)}}}\\
&={{(a-d)^{k(n)}((-1)^{m(n)-k(n)}y^{m(n)}\tilde{\psi}_{n}^{\prime}({{1}\over{y}}))}\over{(2(1-y))^{m(n)}x^{\gamma(n)}}}
\end{align*}

and by theorem \ref{leading}, $(-1)^{m(n)-k(n)}y^{m(n)}\tilde{\psi}_{n}^{\prime}({{1}\over{y}}) \in \mathbb{Z}[a,d,y]$.
 
By the birational equivalence, for any $(x, y) \in E$, \[\psi_{n}(x,y) = 0 \Leftrightarrow \psi_{n}^{\prime}\left(x, {{1}\over{y}}\right)=0\]

so \[\tilde{\psi}_{n}(y) = 0 \Leftrightarrow (-1)^{m(n)-k(n)}y^{m(n)}\tilde{\psi}_{n}^{\prime}({{1}\over{y}}) =0\]

which gives \[\tilde{\psi}_{n}(y) = t(-1)^{m(n)-k(n)}y^{m(n)}\tilde{\psi}_{n}^{\prime}({{1}\over{y}})\]
for some $t$. By comparing leading terms using theorems \ref{leading} and \ref{last}, we get $t = 1$, i.e.,
 
\[\tilde{\psi}_{n}(y) = (-1)^{m(n)-k(n)}y^{m(n)}\tilde{\psi}_{n}^{\prime}({{1}\over{y}}).\]

Now, \[\tilde{\psi}_{n}(a,d,y) = \alpha_{m(n)}(a,d)y^{m(n)} + \alpha_{m(n)-1}(a,d)y^{m(n)-1}+ \dots + \alpha_{1}(a,d)y+\alpha_{0}(a,d)\]
and \[\tilde{\psi}^{\prime}_{n}(a,d,y) = \alpha_{m(n)}(d,a)y^{m(n)} + \alpha_{m(n)-1}(d,a)y^{m(n)-1}+ \dots + \alpha_{1}(d,a)y+\alpha_{0}(d,a).\]

Recall (lemma \ref{hom}) that each of the $\alpha_i$ is homogeneous in $a$ and $d$ of degree $m(n) -k(n)$, so 
\[(-1)^{m(n)-k(n)}\tilde{\psi}^{\prime}_{n}(a,d,y) = \alpha_{m(n)}(-d,-a)y^{m(n)} + \alpha_{m(n)-1}(-d,-a)y^{m(n)-1}+ \dots + \alpha_{1}(-d,-a)y+\alpha_{0}(-d,-a)\]

and \begin{align*}
(-1)^{m(n)-k(n)}y^{m(n)}\tilde{\psi}_{n}^{\prime}({{1}\over{y}}) &= \alpha_{m(n)}(-d,-a) + \alpha_{m(n)-1}(-d,-a)y+ \dots \\
&\quad +\alpha_{1}(-d,-a)y^{m(n)-1}+\alpha_{0}(-d,-a)y^{m(n)}\\
&=\tilde{\psi}_{n}^*(-d,-a,y).
\end{align*}
Hence, $\tilde{\psi}_{n}(a,d,y) = \tilde{\psi}_{n}^*(-d,-a,y).$ 
\hfill $\ \Box$
\bigskip

\section{Another Approach to Division Polynomials}\label{alt}

\subsection{Rephrasing the addition laws}


Let $(x_{+}, y_{+}) = (x_1, y_1)+(x_2, y_2) $, $(x_{-}, y_{-}) = (x_1, y_1)-(x_2, y_2)$

\begin{thm} \label{xadd} \[x_{+}=\frac{x_1y_2(1-dx_2^2)+x_2y_1(1-dx_1^2)}{1-adx_1^2x_2^2}\]
\end{thm}
\emph{Proof:} \begin{align*} x_{+} &= \frac{(x_1y_2+x_2y_1)(1-dx_1x_2y_1y_2)}{1-d^2x_1^2x_2^2y_1^2y_2^2}\\
&=\frac{x_1y_2(1-dx_2^2y_1^2)+x_2y_1(1-dx_1^2y_2^2)}{1-d^2x_1^2x_2^2y_1^2y_2^2}\\
&=\frac{x_1y_2(1-dx_2^2\frac{1-ax_1^2}{1-dx_1^2})+x_2y_1(1-dx_1^2\frac{1-ax_2^2}{1-dx_2^2})}{1-d^2x_1^2x_2^2(\frac{1-ax_1^2}{1-dx_1^2})(\frac{1-ax_2^2}{1-dx_2^2})}\\
&=\frac{(1-d(x_1^2+x_2^2)+adx_1^2x_2^2)(x_1y_2(1-dx_2^2)+x_2y_1(1-dx_1^2))}{(1-dx_1^2)(1-dx_2^2)-d^2x_1^2x_2^2(1-ax_1^2)(1-ax_2^2)}\\
&=\frac{(1-d(x_1^2+x_2^2)+adx_1^2x_2^2)(x_1y_2(1-dx_2^2)+x_2y_1(1-dx_1^2))}{(1-d(x_1^2+x_2^2)+adx_1^2x_2^2)(1-adx_1^2x_2^2)}\\
&=\frac{x_1y_2(1-dx_2^2)+x_2y_1(1-dx_1^2)}{1-adx_1^2x_2^2}
\end{align*} \hfill $\Box$

Notes: If $ad$ is a nonsquare in $K$, it is immediate that the above addition law is \emph{complete} (in the sense of \cite{BL}). It is also straightforward to see that \[x_{-}=\frac{x_1y_2(1-dx_2^2)-x_2y_1(1-dx_1^2)}{1-adx_1^2x_2^2},\] and thus the following theorem holds.

\begin{thm} \label{xplusmin} \[x_{+}+x_{-} = \frac{2x_1y_2(1-dx_2^2)}{1-adx_1^2x_2^2}.\]
\end{thm}

Analogously: \[y_{+} = \frac{(a-d)y_1y_2-(a-dy_1^2)(a-dy_2^2)x_1x_2}{a-d(y_1^2+y_2^2)+dy_1^2y_2^2}\]
\emph{Proof:} \begin{align*} y_{+} &= \frac{(y_1y_2-ax_1x_2)(1+dx_1x_2y_1y_2)}{1-d^2x_1^2x_2^2y_1^2y_2^2}\\
&= \frac{y_1y_2(1-adx_1^2x_2^2)-x_1x_2(a-dy_1^2y_2^2)}{1-d^2x_1^2x_2^2y_1^2y_2^2}\\
&= \frac{y_1y_2((a-dy_1^2)(a-dy_2^2)-ad(1-y_1^2)(1-y_2^2))-x_1x_2(a-dy_1^2y_2^2)(a-dy_1^2)(a-dy_2^2)}{(a-dy_1^2)(a-dy_2^2)-dy_1^2y_2^2(1-y_1^2)(1-y_2^2)}\\
&= \frac{(a-d)(a-dy_1^2y_2^2)y_1y_2-(a-dy_1^2)(a-dy_2^2)(a-dy_1^2y_2^2)x_1x_2}{(a-dy_1^2y_2^2)(a-d(y_1^2+y_2^2)+dy_1^2y_2^2)}\\
&= \frac{(a-d)y_1y_2-(a-dy_1^2)(a-dy_2^2)x_1x_2}{a-d(y_1^2+y_2^2)+dy_1^2y_2^2}
\end{align*} \hfill $\Box$

Thus \[y_{-} = \frac{(a-d)y_1y_2+(a-dy_1^2)(a-dy_2^2)x_1x_2}{a-d(y_1^2+y_2^2)+dy_1^2y_2^2}\] and 

\begin{thm} \label{yplusmin} \[y_{+}+ y_{-} = \frac{2(a-d)y_1y_2}{a-d(y_1^2+y_2^2)+dy_1^2y_2^2}\]
\end{thm}

\subsection{Recursion formulae}
Motivated by the polynomials studied by Abel in proving his theorem on the $n$-division points of the lemniscate \cite{Abel} (and see also Cox \cite{Cox}), we use the above addition formulae to derive a new set of polynomials defined by a recursion to specify the $n$th multiple of a point. 
From here on we denote the $x$-coordinate of $[n](x,y)$ by $x_n$, and the $y$-coordinate by $y_n$.

\begin{thm} \label{PQ} \begin{equation*}x_n =\left\{
\begin{array}{rl}
\frac{xyP_n(x^2)}{Q_n(x^2)} \ &\text{ if $n$ is even}\\
\\
\frac{xP_n(x^2)}{Q_n(x^2)} \ &\text{ if $n$ is odd}
\end{array} \right.
\end{equation*}
where $P_n(t)$, $Q_n(t) \in \mathbb{Z}[t]$ are defined by:
\[P_1(t) = 1, \ Q_1(t) = 1,\quad P_2(t) = 2(1-dt), \ Q_2(t) = 1-adt^2\]

\begin{equation*} P_{n+1}(t) = \left\{
\begin{array}{rl}
&2(1-at)(1-dt)P_nQ_{n-1}Q_{n}-P_{n-1}((1-dt)Q_n^2-adt^2(1-at)P_n^2)\ \text{ if $n$ is even}\\
\\
&2(1-dt)P_nQ_{n-1}Q_{n}-P_{n-1}(Q_n^2-adt^2P_n^2)\quad \text{ if $n$ is odd}
\end{array} \right.
\end{equation*}
\newline{} 

\begin{equation*} Q_{n+1}(t) = \left\{
\begin{array}{rl}
&Q_{n-1}((1-dt)Q_n^2-adt^2(1-at)P_n^2)\ \text{ if $n$ is even}\\
\\
&Q_{n-1}(Q_n^2-adt^2P_n^2)\ \text{ if $n$ is odd}
\end{array} \right.
\end{equation*}

\end{thm}

Note that ($P_{n+1}, \ Q_{n+1}$) is generated by a recursion on ($P_n, \ Q_n$) and ($P_{n-1},\ Q_{n-1}$), as distinct from the recursions on various polynomials of index $\sim \frac{n}{2}$ as in theorem \ref{main}.\\

\emph{Proof:} By induction on $n$. The claim is true for $n = 1$, and, by Theorem \ref{xadd}, for $n = 2$. Assume the claim is true for $n$, $n-1$. Then, by Theorem \ref{xplusmin}, 
\[x_{n+1}+x_{n-1} = \frac{2x_ny(1-dx^2)}{1-adx_n^2x^2}\]

Case 1: $n$ even

\begin{align*}
x_{n+1} &= \frac{2xy^2\frac{P_n}{Q_n}(1-dx^2)}{1-adx^4y^2\frac{P_n^2}{Q_n^2}}-\frac{xP_{n-1}}{Q_{n-1}}\\
&=\frac{2xy^2P_nQ_n(1-dx^2)}{Q_n^2-adx^4y^2P_n^2}-\frac{xP_{n-1}}{Q_{n-1}}\\
&=\frac{2x(1-ax^2)(1-dx^2)P_nQ_n}{(1-dx^2)Q_n^2-adx^4(1-ax^2)P_n^2}-\frac{xP_{n-1}}{Q_{n-1}}\\
&=\frac{x(2(1-ax^2)(1-dx^2)P_nQ_{n-1}Q_{n}-P_{n-1}((1-dx^2)Q_n^2-adx^4(1-ax^2)P_n^2))}{Q_{n-1}((1-dx^2)Q_n^2-adx^4(1-ax^2)P_n^2)}
\end{align*}
proving the claim for the case of $n$ being even.

Case 2: $n$ odd

\begin{align*}
x_{n+1} &= \frac{2xy\frac{P_n}{Q_n}(1-dx^2)}{1-adx^4\frac{P_n^2}{Q_n^2}}-\frac{xyP_{n-1}(x^2)}{Q_{n-1}(x^2)}\\
&=\frac{2xyP_nQ_n(1-dx^2)}{Q_n^2-adx^4P_n^2}-\frac{xyP_{n-1}}{Q_{n-1}}\\
&=\frac{xy(2(1-dx^2)P_nQ_{n-1}Q_{n}-P_{n-1}(Q_n^2-adx^4P_n^2))}{Q_{n-1}(Q_n^2-adx^4P_n^2)}
\end{align*}
Proving the claim for the case of $n$ being odd, and thus, by induction, the theorem. \hfill $\Box$

Equally, one could rephrase the previous theorem as a recursion of rational functions.

\begin{thm} \label{alpha}  \begin{equation*}x_n =\left\{
\begin{array}{rl}
xy\alpha_n(x^2) \ &\text{ if $n$ is even}\\
\\
x\alpha_n(x^2) \ &\text{ if $n$ is odd}
\end{array} \right.
\end{equation*}
where $\alpha_n(t)$ are defined by:
\[\alpha_1(t) = 1,\quad \alpha_2(t) = \frac{2(1-dt)}{1-adt^2},\]

\begin{equation*} \alpha_{n+1}(t) = \left\{
\begin{array}{rl}
\frac{2(1-at)(1-dt)\alpha_n}{(1-dt)-adt^2(1-at)\alpha_n^2}-\alpha_{n-1}\ &\text{ if $n$ is even}\\
\\
\frac{2(1-dt)\alpha_n}{1-adt^2\alpha_n^2}-\alpha_{n-1}\ &\text{ if $n$ is odd}
\end{array} \right.
\end{equation*}
\newline{} 
\end{thm}

\emph{Proof:} Similar \hfill $\Box$

We can also express $x_n$ in terms of $y$, and $y_n$ in terms of $y$ or $x$. For brevity's sake, we omit these formulae.

\subsection{Recovering the $y$ coordinate}
The formulae above can be used to perform $x$-coordinate-only arithmetic (cf Montgomery ladder, \cite{Mont}). For this purpose, we manipulate Theorem \ref{xadd} and the analogous result for $y_{+}$ to get

\begin{thm} \label{recovery}
\[y_n = \frac{x_{n-1}(1-adx^2x_n^2)+x_ny(1-dx^2)}{x(1-dx_n^2)}\]
\[x_n = \frac{y_{n-1}(a-d(y^2+y_n^2)+dy^2y_n^2)-(a-d)yy_n}{(a-dy^2)(a-dy_n^2)}\]
\end{thm}
\emph{Proof:} Immediate from \[x_{+}=\frac{x_1y_2(1-dx_2^2)+x_2y_1(1-dx_1^2)}{1-adx_1^2x_2^2}\] and \[y_{+} = \frac{(a-d)y_1y_2-(a-dy_1^2)(a-dy_2^2)x_1x_2}{a-d(y_1^2+y_2^2)+dy_1^2y_2^2}.\] \hfill $\Box$

\section{Acknowledgement} We thank Dan Bernstein and Tanja Lange for their advice, and for directing us to the 3rd volume of Gauss's \emph{Werke} \cite{Gauss}, discussed in Section \ref{gauss}.


\end{document}